\documentclass[12pt]{iopart}
\usepackage{iopams, setstack}
\usepackage{amssymb}
\usepackage{graphicx}
\usepackage{bm}

\newtheorem{theorem}{Theorem}
\newtheorem{lemma}[theorem]{Lemma}

\newtheorem{definition}{Definition}
\newtheorem{algorithm}{Algorithm}
\def\qed{{\hfill$\Box$}}

\def\Proof{{\noindent \em Proof.}\ }

\begin{document}

\title{Degree-based graph construction}

\author{
Hyunju Kim$^{1}$,
Zolt\'an Toroczkai$^{1,2}$,
P\'eter L. Erd\H{o}s$^{2}$,
Istv\'an Mikl\'os$^{2}$ and
L\'aszl\'o A. Sz\'ekely$^{3}$ 
}

\address{
$^{1}$Interdisciplinary Center for Network Science and Applications (iCeNSA), and Department of Phsyics, University of Notre Dame, Notre Dame, IN, 46556, USA\\
$^{2}$Alfr\'ed R\'enyi Institute of Mathematics, Hungarian Academy of
Sciences, Budapest, PO Box 127, H-1364, Hungary\\
$^{3}$Department of Mathematics, University of South Carolina
Columbia, SC 29208, USA}

\begin{abstract}
Degree-based graph construction is an ubiquitous  problem 
in  network modeling \cite{StrDynNetw_NBW06,
PhysRep_BLMCH06}, ranging from social sciences to chemical compounds and
biochemical reaction networks in the cell.
This problem includes existence, enumeration, exhaustive construction and sampling
questions with aspects that are still open today.  Here we give necessary and sufficient conditions for
a sequence of nonnegative integers to be realized as a simple graph's degree
sequence, such that a given (but otherwise arbitrary) set of 
connections from a arbitrarily given node are avoided.  We then use this result to
present a swap-free algorithm that builds {\em all} simple graphs 
realizing a given degree sequence. In a wider context, we show that our result
provides a greedy construction method to build all the $f$-factor subgraphs
(Tutte, \cite{CanJMath_T52}) embedded
within $K_n\setminus S_k$, where $K_n$ is the complete graph and $S_k$
is a star graph centered on one of the nodes.
\end{abstract}



\maketitle 

\section{Introduction and Summary} \label{intro}

In network modeling of  complex systems \cite{StrDynNetw_NBW06,
PhysRep_BLMCH06}, one usually defines a
graph with components of the system being represented by the nodes, and the 
interactions amongst the components being represented as the links 
(edges) of this graph. This graph is usually inferred from empirical observations
of the system and it is uniquely determined if one can specify all the
connections  in the graph. Occasionally, however, the data 
available from the system is incomplete, and one cannot uniquely determine
this graph. In this case there will be a {\em set} ${\cal G}$
of graphs satisfying the data, and one is faced with the following problems: 1) 
construct a graph from ${\cal G}$; 2) count the number of elements (graphs)
in ${\cal G}$; 3) construct all graphs from ${\cal G}$ and 4) construct a {\em typical} element
of ${\cal G}$, often interpreted as a uniform random sample taken from ${\cal G}$. 
Problems 1),3) and 4) are construction type problems, whereas 2) is an
enumeration type problem. 
In this paper we restrict ourselves to simple, undirected
graphs, that is, any edge connects a single pair of distinct nodes (no hypergraph, no
self-loops) and there is at most one edge incident on any pair of nodes (no parallel 
or multiple edges). 

A rather important and typical situation is when the empirical data  specifies only the 
degrees of the nodes, in form of a sequence $\bm{d}=\{d_1,d_2,\ldots,d_n\}$ of positive
integers, $d_i \geq 1$, $1 \leq i \leq n$. (We exclude zero degree nodes, 
since they represent isolated points.) In the following we will refer to 
such cases as ``degree based" construction (enumeration)
problems.
There are numerous examples, we will mention only a few here.
In epidemics studies of sexually transmitted diseases \cite{Nature_LEASA01} the 
data  collected is from anonymous surveys, where the 
individuals specify the {\em number} of different
partners they have had in a given period of time, without revealing their identity. 
In this case, the epidemiologist is faced with constructing the most typical contact graph obeying 
the empirical degree sequence. Another example comes from chemistry where the task is
to determine the total number of structural isomers of chemical compounds, such as alkanes. 
In this case nodes represent chemical elements (atoms) in the compound and a link represents 
a chemical bond. In the case of alkanes the bond can be interpreted as a single link in the corresponding
graph (no double bonds).
Since the valence of an atom is fixed, the formula of an alkane
such as $C_4H_{10}$ (butane) will specify only the degree sequence. Knowing all the possible 
graphs with this degree sequence provides a starting point from which the  feasible
structures can be inferred. In particular, butane has 2,
octane ($C_8H_{18}$) has 18, $C_{20}H_{42}$ has 366,319 isomers, etc. 
Degree-based graph construction is also found in many other network modeling problems,
such as communications (Internet, www, peer-to-peer networks,) biology (metabolic
networks, gene transcription, etc.) ecology (food webs) and social networks.  

It is easy to see that not all integer sequences can be realized as the degrees 
of a simple graph (the existence problem). 
For example, while $\{2,1,1\}$ and $\{2,2,2\}$ are the degree sequences of a 
path ($P_3$, \textbullet--\textbullet--\textbullet  ) and a triangle, there 
is no simple graph with degree sequence $\{3,2,1\}$ or $\{1,1,1\}$ or $\{4,4,2,1,1\}$.  
Let $G(V,E)$ denote a simple graph where $V = \{v_1,v_2,\ldots,v_n \}$ denotes the set of nodes
and $E$ the set of edges. 
 Consider a  sequence of positive integers $\bm{d}=\{d_1,d_2,\ldots,d_n\}$ 
arranged in decreasing order, $d_1 \geq d_2 \geq \ldots \geq d_n$ 
(for convenience reasons, only). If there is a simple graph $G(V,E)$ with degree sequence $\bm{d}$, 
then we call the sequence $\bm{d}$ a {\em graphical sequence} 
and in this case we also say that $G$ {\em realizes}  $\bm{d}$. A second observation is that 
given a graphical $\bm{d}$ (and thus, we know that a simple graph $G$ exists with this degree sequence), 
careless connections of pairs of nodes may not result in a  simple graph. For example, 
consider the sequence $\{2,2,2,2\}$  which is graphical (4-cycle). Making the
connections $\{(v_1,v_2),(v_1,v_3),(v_2,v_3)\}$ however, will force us to make a self-loop 
$\{(v_4,v_4)\}$. The degree-based graph construction problem for 
{\em simple undirected labelled graphs} thus can be 
announced as follows:

\medskip\noindent {\bf Degree-based Graph Construction:} \\
Given  a sequence of integers $\bm{d}=\{d_1,d_2,\ldots,d_n\}$, $d_1 \geq \ldots \geq d_n \geq 1$,
\begin{itemize}
\item[A)] Is there a simple graph $G(V,E)$ on $n$-nodes realizing $\bm{d}$?
\item[B)] If the answer to A) is yes, how can we build such a graph?
\item[C)] Can we build {\em all} such graphs?
\item[D)] Let ${\cal G}(\bm{d})$ be the set of all such graphs. 
How can we sample at uniform from  ${\cal G}(\bm{d})$?
\end{itemize}

There are two well-known theorems that answer 
question A) above, namely the Erd\H{o}s-Gallai theorem \cite{EG60} and the Havel-Hakimi theorem
\cite{H55,H62},
the latter also giving a construction algorithm for a graph with degree sequence $\bm{d}$
and thus answering question B) as well (see Section \ref{prev}). 
In principle, problem  C) can be resolved 
via the method of edge swaps starting from the graph produced by the Havel-Hakimi procedure  from B) 
(called an HH-graph from now on) and book-keeping the swaps (which gets rather involved).  Given two edges
$(v_1,v_2)$ and $(v_3,v_4)$, they can be swapped into $(v_1,v_3)$ and $(v_2,v_4)$, or
$(v_1,v_4)$ and $(v_2,v_3)$ leaving the degree sequence unchanged. Due to a theorem by
Ryser \cite{CanJMath_R57}, if $G_1$ and $G_2$ are two simple graphs with identical degree
sequences, then there is a sequence of edge swaps that
transforms one into another \cite{LinAlgAppl_B80, SIAM_JAlgDiscMeth_T82}. 
Edge swapping is also at the basis of all sampling algorithms attempting
to answer D), using a Markov Chain Monte Carlo approach, the literature of which is too extensive
to be reviewed here. The basic idea is to keep swapping edges until the memory of the 
initial condition (HH-graph) is lost and one produces a (pseudo)-random instance. This sampling method is approximative and it is not well controlled in general (except for some specific sequences). A simple and direct (swap-free)
construction method to produce a uniformly sampled random graph from ${\cal G}
(\bm{d})$ was presented by Molloy and Reed (M-R) \cite{RandStrucAlg_MR95} 
(see Section \ref{discussion})  and 
subsequently used to generate graphs with given 
degree sequences \cite{PhysRevE_NSW01}, including 
those described by power-law degree distributions 
\cite{Proc32ACMSympTheorComp_ACL00}. 
The problem with the M-R algorithm is that it can 
become very slow due to rejections caused 
by self-loops and parallel edges (see Section \ref{discussion}).

Here we present a  new approach to degree-based graph construction. First we prove our main 
result that gives the sufficient and necessary conditions for a sequence of integers to be 
graphical such 
that a given (but otherwise arbitrary) set of  connections from a given (but otherwise arbitrary) 
node are avoided. We then show how to use this result to present an algorithm that builds 
{\em all graphs} from ${\cal G}(\bm{d})$ (question C)).  It is important to emphasize, that our 
algorithm does not use edge swaps, it is a direct construction method. Lastly, (Section \ref{discussion}) 
we show how our 
result  improves on the M-R method of uniform sampling, allowing to reject some of
the samples without actually getting to the point where the multi-edge 
conflicts would actually occur (see Section \ref{discussion}). 
We also make a connection with the renowned Tutte's $f$-factor theorem 
\cite{CanJMath_T52,CanJMath_T54}, showing that our result provides a greedy algorithm for constructing all
$f$-factors in case of $K_n\setminus S_k$ where $K_n$ is the complete graph on $n$ nodes and
$S_k$ is a star graph with $k$ leaves centered on some arbitrary node.

The paper is organized as follows: Section \ref{prev} recalls known fundamental theorems
for graph construction, and via a simple counter-example it shows that HH
is not sufficient to build all graphs from ${\cal G}(\bm{d})$; Section \ref{main} 
presents our main theorem
with its proof and Section \ref{all} describes the algorithm for building all graphs 
in ${\cal G}(\bm{d})$; Section \ref{discussion} is devoted to discussions.

\section{Previous results} \label{prev}

For simplicity of the notation
we will identify node $v_i$ by the integer $i$. 
There are two well-known necessary and sufficient conditions for a
sequence of nonnegative integers to be graphical: one was given
independently by Havel \cite{H55} and Hakimi \cite{H62} while the
other is due to Erd\H{o}s and Gallai \cite{EG60}. We now announce 
these results for later reference, however, without proof, those
can be found in the corresponding references.
\begin{theorem}[Hakimi-Havel, HH]\label{tm:HH}
There exists a simple graph with degree sequence $d_1 > 0$, $d_2 \ge ... \ge
d_n > 0$ if and only if there exists one with degree
sequence $d_2-1,\ldots,d_{d_1+1}-1,d_{d_1+2},\ldots,d_n$.
\end{theorem}
\begin{theorem}[Erd\H os-Gallai, EG]\label{tm:EG}
Let $ d_1 \ge d_2 \ge ...... \ge d_n>0$ be integers. Then they are
the degree sequence of a simple graph if and only if \\
{\rm (i)} $d_1+...+d_n$ is even \\
{\rm (ii)} for all $k=1,...,n-1$ we have
\begin{equation}
\sum_{i=1}^k d_i \le k(k-1) + \sum_{i=k+1}^n \min\{k, d_i\}\;. \label{ineq}
\end{equation}
\end{theorem}
Note that Theorem~\ref{tm:HH} provides a greedy algorithm (the HH-algorithm) 
to generate an
actual graph with the given degree sequence $\bm{d}$ while
Theorem~\ref{tm:EG} is an existence result. Tripathi and Vijay have
recently shown \cite{DiscrMath_TV03} that it is enough to check the inequalities (\ref{ineq})
for $1\leq k \leq s$, where $s$ is determined by $d_s \geq s$,
$d_{s+1} < s+1$, that is only as many times as many distinct terms are in the
degree sequence.

In the following we will imagine the given degree sequence as a
collection of {\em stubs}: at each vertex $i$ there are $d_i$ stubs (``half-edges"),
anchored at the vertex, but the other ends are free. Connecting two
stubs at two distinct nodes will form an edge between those nodes.
We will call {\em the residual degree} the number of current stubs of a node.

The HH-algorithm for constructing a graph realizing a graphical sequence 
$\bm{d}$ proceeds as follows: connect all stubs of a node  to
nodes that have the largest residual degrees and repeat until
no stubs are left. It is important to emphasize, that one can choose any node
to connect its stubs, as long as we connect all its stubs to the other nodes with
the largest residual degrees.  Clearly, if we always choose a high degree node (from the 
residual sequence) to connect its stubs, the HH-algorithm will create
a graph in which high degree nodes tend to be connected to other
high degree nodes, called assortative property \cite{PhysRevLett_M02}. 
However, if we always pick a node with low (residual) degree to connect,
we will likely obtain a graph with dissassortative property  \cite{PhysRevE_M03}. 
The HH-theorem is also a consequence (as a corollary) of our main 
result, see Section \ref{main}.
Nevertheless, this is still not enough 
to produce {\em all} graphs realizing a graphical sequence!
To see that, consider the graphical sequence $\bm{d} = \{3,3,2,2,2,2,2,2\}.$

\begin{figure}[htbp]
\bigskip \centering \vspace*{-0.5cm}
\includegraphics[width=2in]{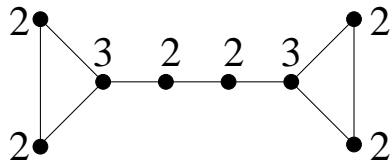}
\caption{\em This graph cannot be obtained by the 
Havel-Hakimi procedure. The integers indicate node
degrees.}\label{fig:nonghh}
\end{figure}

If the first node to connect is a node with degree 3, then the
HH-algorithm connects it to the other node with degree 3 (highest
degree). If the first node to connect has degree 2, then
the HH-algorithm connects both its stubs to nodes with degree 3.
However, the graph in Fig. \ref{fig:nonghh} does not have any of the 
connections just mentioned (a $3 - 3$ or $3 - 2 - 3$ connection), and thus
it cannot be constructed with the HH-algorithm. In the next section
we introduce a theorem that allows us to construct {\em all} labelled graphs
with a given degree sequence.

The above results  are naturally placed in the larger context of 
Tutte's famous $f$-factor theorem \cite{CanJMath_T52}.  Given an integer
function $f:V\rightarrow \mathbb{N}\cup \{0\}$, 
the {\em $f$-factor} of a given simple graph $G(V,E)$ is a subgraph 
 $H$ of $G$ such that $d_H(v)=f(v)$ for all $v\in V$. Here $d_H(v)$ is
 the degree of $v$ within $H$.  Tutte gave sufficient and necessary  
 conditions for the existence 
of an $f$-factor for $G$  \cite{CanJMath_T52}, and later connected this to the problem of finding
perfects matchings in bipartite graphs \cite{CanJMath_T54}. 
It is not hard to see, that taking $G=K_n$, that is the complete graph on
$n$-nodes, the $f$-factor problem is exactly question A) of the degree-based
construction problem with $\bm{d} = \{f(v_1),\ldots,f(v_n)\}$. 
In this sense, the HH-algorithm is a greedy 
method for constructing an $f$-factor on $K_n$.
 
\section{Graphical sequences with constraints from a single node} \label{main}

Before we can announce and prove our main result, we need to
introduce a number of definitions and observations. 
\begin{definition}
Let $A(i)$ be an increasingly ordered set of $d_i$ distinct nodes
associated with node $i$: $A(i)= \{a_k \; | \; a_k \in V,\;\; a_k
\ne i, \;\;\forall k,~ 1\le k \le d_i \}$.
\end{definition}
Usually, this set will represent the set of nodes adjacent to node
$i$ in some graph $G$, therefore we will refer to $A(i)$ as an {\em
adjacency set} of $i$.
\begin{definition}\label{2sets}
If for two adjacency sets $A(i) = \{\ldots,a_k,\ldots\}$ and $B(i) =
\{\ldots,b_k,\ldots\}$ we have $b_k \le a_k$ for all $1 \le k \le
d_i$, we say that $B(i) \le A(i)$.
\end{definition}
In this case we also say that $B(i)$ is ``to the left" of $A(i)$.
\begin{definition} \label{reduced}
Let $d_1\ge d_2 \ge\dots\ge d_n \ge 1 $ be a graphical sequence, and
let $A(i)$ be an adjacency set of node $i$. The degree sequence
reduced by $A(i)$ is defined as:
\begin{eqnarray}
d'_k\big|_{A(i)} =\left\{\begin{array}{lll} d_k - 1 & \quad \mbox{if }
       & k \in A(i)\\
       d_k & \quad \mbox{if } & k \in [1,n] \setminus (A(i) \cup \{i\}) \\
       0 & \quad \mbox{if} & k = i\; .
   \end{array}\right. \label{red}
\end{eqnarray} 
\end{definition}
\noindent In other words, if $A(i)$ is the set of adjacent nodes to $i$ in the graph
$G$, then the reduced degree sequence $\bm{d'}\big|_{A(i)}$ is
obtained after removing node $i$ with all its edges from $G$.
\begin{lemma} \label{DS}
Let $\{d_1,\ldots,d_j,\ldots,d_k,\ldots, d_n\}$ be a non-increasing
graphical sequence and assume $d_j > d_k.$ Then the sequence
$\{d_1,\ldots,d_j-1,\ldots,d_k+1,\ldots, d_n\}$ is also graphical
(not necessarily ordered).
\end{lemma}
\Proof Since $d_j > d_k$, there exists a node $m$ connected to node
$j$, but not connected to node $k$. Let us cut edge $(m,j)$ and
remove the disconnected stub of $j$. If we add one more stub to $k$,
and connect this new stub to the disconnected stub of $m$, then we
can see that the new graph is also simple with degree sequence $\{
d_1,~ d_2,~\dots~,~ d_j - 1,~\dots~,~ d_k + 1 ,~\dots~,~d_{n} \}$.
\hfill \ \qed
\begin{lemma}\label{ShS}
Let $\bm{d}=\{d_1,d_2,\ldots, d_n\}$, be a non-increasing graphical
sequence, and let $A(i)$, $B(i)$ be two adjacency sets for some node $i
\in V$, such that $B(i) \le A(i)$. If
the degree sequence  reduced by $A(i)$ (that is $\bm{d'}\big|_{A(i)}$) is
graphical, then the degree sequence 
reduced by $B(i)$ (that is $\bm{d'}\big|_{B(i)}$) is also graphical.
\end{lemma}
\Proof Let $A(i)=\{\ldots, a_k, \ldots \}$ and
$B(i)=\{\ldots,b_k,\ldots\}$, $k=1,\ldots,d_i$. Consider the
adjacency set $B^1(i) = \{b_1,a_2,a_3,\ldots,a_{d_i}\}$ (we replaced
node $a_1$ by node $b_1 \le a_1$). If $b_1 = a_1$ then there is
nothing to do, we move on (see below). If $b_1 < a_1$ then
conditions in Lemma \ref{DS} are fulfilled. Namely, $b_1 < a_1$
implies $d_{b_1} \geq d_{a_1} > d_{a_1}-1$ and we know that the
sequence $\bm{d'}\big|_{A(i)} = \{\ldots,d_{b_1},\ldots,
d_{a_1}-1,\ldots,d_{a_2}-1, \ldots\}$ is graphical by assumption.
Thus, according to Lemma \ref{DS}, the sequence
$\{\ldots,d_{b_1}-1,\ldots, d_{a_1},\ldots,d_{a_2}-1, \ldots\}$
is also graphical, that is
the one reduced by the set $B^1(i)$. Next, we will proceed by
induction. Consider the adjacency set
$B^m(i)=\{b_1,\ldots,b_m,a_{m+1},a_{m+2},\ldots,a_{d_i}\}$ and
assume that the degree sequence reduced by it (from $\bm{d}$)
is graphical. Now, consider the adjacency set
$B^{m+1}(i)=\{b_1,\ldots,b_{m+1},a_{m+2},a_{m+3},\ldots,a_{d_i}\}$
(replaced $a_{m+1}$ by $b_{m+1}$). If $b_{m+1} < a_{m+1}$, Lemma
\ref{DS} can be applied again since $b_{m+1} < a_{m+1}$ implies
$d_{b_{m+1}} \geq d_{a_{m+1}} > d_{a_{m+1}}-1$, showing that the
sequence reduced by $B^{m+1}(i)$ is also graphical. The last
substitution ($m+1 = d_i$) finishes the proof. \hfill \ \qed
\begin{definition}\label{Lset} 
Let $\bm{d}=\{d_1,d_2,\ldots, d_n\}$ be a decreasing graphical
sequence and consider an arbitrary node $i \in V$, and
 an arbitrarily fixed integer $m$ with $0 \leq m \le
n-1-d_i$. Let us fix a set of nodes $X(i) =
\{j_1,\ldots,j_m\} \subset V\setminus \{i\}$ and consider the set
$L(i)=\{l_1,\ldots,l_{d_i}\}$ containing the $d_i$ lowest index nodes
not in $X(i)$ and different from $i$. We call $L(i)$ the leftmost
adjacency set of $i$ restricted by $X(i)$. Accordingly, we call the
set of nodes $X(i)$ the set of forbidden connections for $i$.
\end{definition}
\begin{lemma}\label{LY}
If $\bm{d}=\{d_1,d_2,\ldots, d_n\}$ is a decreasing graphical
sequence, and $Y(i) = \{y_1,\ldots y_{d_i}\}$ is an adjacency set
disjoint from $X(i) \cup \{i\}$, then $L(i) \leq Y(i)$.
\end{lemma}
\Proof This is immediate, since by Definition \ref{Lset},
$l_j \leq y_j$, for all $j=\{1,\ldots,d_i\}$. \hfill \ \qed \\
We are now ready for the main theorem:

\begin{theorem}[Star-constrained graphical sequences] \label{CG}
Let $d_1 \ge d_2 \ge \ldots d_n \geq 1$ be a sequence of integers.
For an arbitrary node $i \in V$ define a set $X(i) =
\{j_1,\ldots,j_m\} \subset V\setminus \{i\}$ with $m \le n-1-d_i$,
and consider $L(i)$, the leftmost adjacency set of $i$ restricted by
$X(i)$. Then the degree sequence $\bm{d} = \{d_1,\ldots,d_n\}$
can be realized by a simple graph $G(V,E)$ in which $(i,j) \not\in
E$, for all $j\in X(i)$, if and only if the degree sequence reduced
by $L(i)$ is graphical.
\end{theorem}
\Proof ``$\Longleftarrow$'' is straightforward: add node $i$ to the
reduced set of nodes, then connect it with edges to the nodes of $L(i)$. Thus
we obtained a graphical realization of $\bm{d}$ in which
there are no connections between $i$ and any node in $X$.\\
``$\Longrightarrow$'' In this case $\bm{d}$ is graphical with no
links between $i$ and $X(i)$, and we have to show that the sequence
obtained from $\bm{d}$ by reduction via $L(i)$ is also
graphical. However, $\bm{d}$ graphical means that there is an
adjacency set $A(i)$ (with $A(i) \cap X(i) = \emptyset $) containing
all the nodes that  $i$ is connected to in $G$. Thus, according
to Lemma \ref{LY}, we must have $L(i) \leq A(i)$. Then, by Lemma
\ref{ShS}, the sequence reduced by $L(i)$ is graphical. \qed 

Note, that, the forbidden set of connections form a star graph $S_{|X(i)|}$
centered on node $i$. 
Also note that considering the empty set as the set of forbidden nodes, 
$X(i) = \emptyset$, we obtain the Havel-Hakimi Theorem \ref{tm:HH} as  
corollary.

\medskip \noindent Informally, Theorem \ref{CG} can be announced as follows:

\medskip\noindent{\em Let $\bm{d}= \{d_1, d_2,$ $\ldots, d_n\}$,
be a decreasing graphical sequence and
let $i$ be a fixed, but arbitrary vertex. Assume we are given
a set of forbidden connections in $V$ incident on  $i$. 
Then there exists a realization of the degree sequence avoiding all 
forbidden connections if and only if
there also exists a realization where $i$ is connected with vertices
of highest degree among the non-forbidden ones.}

\medskip\noindent Since the forbidden connections emanating from a node $i$ form 
a star graph $S_k$, $k= |X(i)|$, Theorem \ref{CG} provides sufficient
and necessary conditions for the existence of an $f$-factor for 
$G = K_n\setminus S_k$. More importantly, it gives a 
{\em greedy algorithm} for finding such an $f$-factor.

\section{Building all graphs from ${\cal G}(\bm{d})$} \label{all}

As we show next, Theorem \ref{CG} provides us with a procedure that allows for the
construction of {\em all} graphs realizing the same degree sequence.

Consider a graphical degree sequence $\bm{d}$ on $n$ nodes.
Certainly, we can produce all graphs realizing this sequence by
connecting all the stubs of a chosen node first, before moving on to
another node with stubs to connect (that is we finish with a node, before
moving on).  In this vein, now choose a node $i$  and 
connect one of its stubs to some
other node $j_1$. Is the remaining degree sequence $\bm{d'} =
\{d_1,\ldots, d_i-1,\ldots ,d_{j_1}-1,\ldots,d_n\}$ still graphical
such that nodes $i$ and $j_1$ avoid another connection in subsequent 
connections of the other stubs? Certainly, as a necessary condition,
$\bm{d'}$ has to be graphical as a sequence, since all subgraphs of
a simple graph are simple, and thus if $G$ is a simple graph realizing 
$\bm{d}$ with  $(i,j_1) \in E$, then after removing this edge, the remaining
graph is still simple. However, it is not sufficient that after making some connections
from a node, the residual sequence to be graphical. One might still be forced
to make multiple edges, as illustrated by the following example. Consider
the graphical sequence $\{2,2,1,1\}$ (the path $P_4$, \textbullet--\textbullet--\textbullet--\textbullet )
as the degrees of the set of nodes $V=\{u,v,x,y\}$ ($d_u=d_v=2$, $d_x=d_y=1$).
Connect nodes $u$ and $v$. 
We certainly have not broken the graphical character yet, since we could still finish
the path by connecting next $u$ to $x$ (or to $y$) and $v$ to $y$ (or to $x$).
The remaining sequence $\{1,1,1,1\}$ as a {\em sequence of
integers} is graphical (two edges). Next, connect node $x$ to node $y$. 
The remaining (residual) sequence is $\{1,1\}$ (emanating from node $u$ and $v$, respectively), 
{\em graphical on its own}, however,
we can no longer connect nodes $u$ and $v$, because the very first connection is
already there. Thus, after making one, or more connections from a node $i$,
how can we check that the next connection from $i$ 
will not break the graphical character?

Theorem \ref{CG} answers this question if we think of the connections already made
from node $i$ as forbidden connections. That is, after the first connection of $i$ to $j_1$
take $\bm{d'}$ as $\bm{d}$ in Theorem \ref{CG} and $X(i)=\{j_1\}$. 
Then, to test whether the sequence reduced by the
corresponding $L(i)$ is graphical we can employ for example the Erd\H{o}s-Gallai
Theorem \ref{tm:EG}, checking all the inequalities, or the 
Havel-Hakimi Theorem \ref{tm:HH}.  If the test fails on the
reduced sequence, one must disconnect $i$ from $j_1$ and reconnect
it somewhere else. The graphical character of the original sequence guarantees that there
is always a $j_1$ where the test will not fail.
If, however, the remaining degree sequence is
graphical with the constraint imposed by $X(i)$, we connect another
stub of $i$ to some other node $j_2$ (different from $j_1$), adding an element
to the forbidden set of connections $X(i)$.  To check whether after the second connection the
remaining sequence is still graphical with the constraint imposed by
the new set $X(i) = \{j_1,j_2\}$, we proceed in exactly the same way,
using Theorem \ref{CG}, repeating the procedure until all the stubs of node $i$
are connected away into edges. After this we can move on to some other
node $i$ (arbitrary) from the remaining set of nodes and repeat the procedure. 
Note that this procedure is not a real procedure in the
sense that it does not prescribe which stubs to connect. It only tells us whether the connection
we just made (by whatever process) has broken the graphical character. Since every
element from ${\cal G}(\bm{d})$ can be realized by some sequence of connections, it is clear
that if we specify a systematic way of going through all the possible connections while employing
Theorem \ref{CG}, we will realize all elements of ${\cal G}(\bm{d})$. However, taking all possible
connections would be very inefficient. Next we present a 
version of a more economical algorithm that constructs every labelled graph with 
degree sequence $\bm{d}$, and only once. For simplicity of the notation
we will call the test for the graphical character via 
Theorem \ref{CG}, the ``CG test'' (constrained graphicality test). 
The algorithm also exploits Lemma \ref{ShS}, which  guarantees preservation of graphicality for all adjacency sets {\em to the left} of a graphical one, thus avoiding costlier checks with EG or HH theorems for those adjacency sets. Clearly, a labeled graph
can be characterized by the sequence of its adjacency sets $G = \{A(1),\ldots,A(n)\}$. This algorithm
creates all the possible adjacency sets for node 1, then {\em for each one of those} repeats the same
procedure on the reduced sequence by that adjacency set (in sense of Definition \ref{reduced})
of {\em at most $n-1$ nodes}.
\begin{algorithm}[All graphs] \label{uno} Given a graphical sequence 
$d_1 \ge d_2 \ge \ldots d_n $ $ \geq 1$,  
\begin{itemize}
\item[I.] Create the rightmost adjacency set $A_R(1)$ for node 1: 
Connect node $1$ to $n$ (this never breaks graphicality). Let $k=n-1$.
\begin{itemize}
\item[I.1] Connect another stub of $1$  to $k$. Run the CG test. 
\item[I.2] If it fails, make $k = k-1$. Repeat I.1
\item[I.3] If passes, keep (save) the connection, make $k = k-1$, and if $i$ has stubs left, repeat from I.1.
 \end{itemize}
  \item[II.] Create the set ${\cal A}(\bm{d})$ of all adjacency sets 
  of node $1$ that are colexicographically smaller than $A_R(1)$ and  preserve graphicality:
  $${\cal A}(\bm{d}) = 
  \left\{A(1)=\{a_1,\ldots, a_{d_1}\}, a_i \in V \Big| A(1) <_{CL} A_R(1),\;\;
 \bm{d'}|_{A(1)}-\mbox{graphical}  \right\}.$$
 \item[III.] For every $A(1) \in {\cal A}(\bm{d})$ create all graphs from the corresponding 
 ${\cal G}(\bm{d'}|_{A(1)})$ using this Algorithm, where $\bm{d'}|_{A(1)}$ is 
 the sequence reduced by $A(1)$.
\end{itemize}
\end{algorithm}
For simplicity of the notation, we will drop the $(1)$ from $A(1)$, tacitly assuming that it refers to the leftmost node 1. Observe that the ordering relation ``$<$'' in Definition \ref{2sets} is a partial order, while the colexicographic order ``$<_{CL}$'' is a total order over all adjacency sets, however, ``$<$'' implies ``$<_{CL}$''.  It is not hard to see, that $A_R$ is colexicographically the largest (``rightmost'') sequence which still preserves graphicality. 
When constructing ${\cal A}(\bm{d})$, checking graphicality with the EG or HH theorems is only needed for those adjacency sets,  which are incomparable by the ``$<$'' relationship to any of the current elements of ${\cal A}(\bm{d})$, while for the rest graphicality is guaranteed by
Lemma \ref{ShS}.  

\section{Discussion and Outlook} \label{discussion}

Algorithm \ref{uno} proceeds by attempting to connect all stubs of the largest
degree node as much to the right as possible. Depending on the degree sequence, it might happen that the CG test fails many times at step I., until it finds $A_R$. 
However, in that case, $A_R$
is located more towards the higher degree nodes (towards left) and thus the number of adjacency sets that 
preserve graphicality, namely $|{\cal A}|$ is smaller and accordingly, the algorithm has
fewer cases to run through in subsequent  steps. The more  heterogeneous is a degree sequence, the more
likely this will happen. Of course, it only makes sense to produce all graphs from ${\cal G}(\bm{d})$ for small graphs
(chemistry), or graphical sequences that do not admit too many solutions.  An interesting question would then 
be finding the conditions on the sequence $\bm{d}$ that would guarantee {\em a given} upper bound $C$ on
the size of ${\cal G}(\bm{d})$, $\left| {\cal G}(\bm{d})\right| \leq C$.  A possible starting point in this
direction could be Koren's \cite{Koren} characterization of sequences uniquely realizable by a simple graph. 
Sequences that admit only a small number of realizations (by labeled simple graphs) would likely
be ``close" in some sense to these special sequences.

Algorithm \ref{uno} also provides a way to computationally enumerate all the labelled graphs $ \left| {\cal G}(\bm{d})\right| $ 
realizing a degree sequence $\bm{d}$ (problem B) of Section \ref{intro}).  Naturally, the following recursion holds:
$
 \left| {\cal G}(\bm{d})\right| = \sum_{A\in{\cal A}(\bm{d})}  \left| {\cal G}(\bm{d'}|_A)\right|\;.
$
Our graph construction process can be thought of as happening along the branches of a  tree 
${\cal T}(\bm{d})$ of depth at most $n-1$: the internal nodes of this tree on the $k$th level are all the allowed  adjacency sets (from the corresponding ${\cal A}$ set) of the node with the largest residual 
degree (the leftmost node). The reason this tree is of depth at most $n-1$ is because some other nodes (other than the one with the
largest residual degree) in the process might loose all their stubs. A directed path towards a leaf of this tree corresponds to a graphical realization of $\bm{d}$, because
we end up specifying all the adjacency sets along this path.
Based on this, during the realization of the graph, if 
we choose uniformly at random at every level of the tree within the children of a node, from the corresponding 
set ${\cal A}$, the probability of a final realization 
$G$ in this process will be given by the product:
\begin{equation}
P(G) = \prod_{k=0} \left| {\cal A}(\bm{d'}|_{A^{(k)}})\right|^{-1} \label{w}
\end{equation}
where $A^{(k)}$ is the randomly chosen adjacency set of the node with the largest residual degree, 
on the $k$-th level of ${\cal T}(\bm{d})$. By convention $A^{(0)} = \emptyset$, and 
$\bm{d'}|_{A^{(0)}} =\bm{d}$. It is not hard to convince ourselves 
that the distribution in (\ref{w}) is not uniform, and thus, this algorithm cannot be used in this form
to produce uniform samples from ${\cal G}(\bm{d})$. However, Theorem \ref{CG} can be used to improve
on a well-known, direct uniform sampling process, the Molloy-Reed (M-R) algorithm 
\cite{RandStrucAlg_MR95}. In this process, one chooses among all the stubs uniformly at random, irrespectively
to what node they belong to. This is repeated until either a self-loop, or a double edge is created, or a
simple graph is finished. When there is a self-loop or double edge, the process is stopped, and the algorithm starts from
the very beginning.  The CG test can be used along the way to test for graphicality after every connection
just made. In particular, assume that we just connected node $i$ with node $j$. We then run the CG test centered
on node $i$ (adding the $(i,j)$ connection to the forbidden set from $i$). If passes, we run the CG test on node 
$j$ as well. If it fails either on $i$ or $j$, we can stop the process, before actually running into a conflict
later. Running into conflict usually happens towards the end, and this test can save us from
possibly many unnecessary
calls to the random number generator, 
speeding up the sampling algorithm. It is important to note, that if the CG test
passes on both $i$ and $j$, we actually do not know whether graphicality is broken or still preserved at that
stage! 
If the CG test fails, however, we know that graphicality was broken. 
The reason is because Theorem \ref{CG} gives us the sufficient and necessary conditions for the sequence
to be graphical such that no multiple edges will be made with the already existing connections {\em emanating
from the same node}. It gives no such information for connections already made elsewhere!  Of course, for
degree sequences $\bm{d}$ for which there is only a small number of labelled graphs realizing it
(compared to the total number $\prod_{i}d_i!$ of graphs that it produces), the M-R algorithm
would search for needles in a haystack, and other (MCMC) methods will be necessary.

In summary, we have given necessary and sufficient conditions for a sequence of integers $d_1 \geq \ldots
d_n \geq 1$
to be graphical (realizable by simple,  undirected graphs) avoiding multiple edges with an arbitrary star
graph (the forbidden graph) $S_{k}(j)$, $0 \leq k < n-d_j$, centered on a node $j$. In a more general context, our result gives  for the first time a greedy construction for Tutte's $f$-factor subgraphs within 
$K_n \setminus S_{k}(j)$. It would
be desirable if such a {\em greedy} construction existed for an arbitrary forbidden graph $F$, not just for star-graphs,
however, at this point this still seems to be a rather difficult problem. Such an algorithm (if greedy) would 
further speed up the M-R sampling, because it would induce early rejections, as soon as they are made.
Our main theorem also led to a direct and systematic construction algorithm that builds {\em all graphs} 
realizing a given degree sequence $\bm{d}$.   And finally, we mention that these studies 
can be extended to simple directed graphs as well (there are two degree sequences in this case),
however, the computations are considerably more involved, and they will be presented 
separately.

\subsection*{Acknowledgements}

This project was supported in part by  the NSF BCS-0826958 (HK and ZT), HDTRA
201473-35045 (ZT) and by Hungarian Bioinformatics MTKD-CT-2006-042794,  Marie Curie Host Fellowships for Transfer of Knowledge (LAS and ZT). 
ELP was partly supported by OTKA (Hungarian NSF), under contract Nos.
NK62321, AT048826 and  K 68262, and LAS by NSF DMS-0701111. IM was supported by a Bolyai postdoctoral stipend and OTKA grant  F61730.

\end{document}